\documentclass[11pt,reqno]{amsart}
\usepackage{amsmath,amssymb,amsthm,amsfonts,enumerate,hyperref, caption}
\theoremstyle{plain}
\newtheorem{theorem}{Theorem}[section]
\newtheorem{proposition}[theorem]{Proposition}
\newtheorem{cor}[theorem]{Corollary}
\newtheorem{lemma}[theorem]{Lemma}
\newtheorem{conj}[theorem]{Conjecture}

\theoremstyle{definition}

\theoremstyle{remark}

\numberwithin{equation}{section}

\newcommand{\Ric}{\mathrm{Ric}}

\begin{document}
\title{Numerical approximations to extremal toric K\"ahler metrics with arbitrary K\"ahler class}
\author{Stuart James Hall}
\address{Department of Applied Computing, University of Buckingham, Hunter St., Buckingham, MK18 1G, U.K.} 
\email{stuart.hall@buckingham.ac.uk}

\author{Thomas Murphy }
\address{Department of Mathematics, California State University Fullerton, 800 N. State College Bld., Fullerton, CA 92831, USA.}
\email{tmurphy@fullerton.edu}

\maketitle
\begin{abstract}
We develop new algorithms for approximating extremal toric K\"ahler metrics. We focus on an extremal metric on $\mathbb{CP}^{2}\sharp2\overline{\mathbb{CP}}^{2}$, which is conformal to an Einstein metric (the Chen-LeBrun-Weber metric). We compare our approximation to one given by Bunch and Donaldson and compute various geometric quantities.  In particular, we demonstrate a small eigenvalue of the scalar Laplacian of the Einstein metric which gives a numerical evidence that the Einstein metric is conformally unstable under the Ricci flow. 
\end{abstract}
\section{Introduction}
This article develops new methods for numerically approximating extremal toric K\"ahler metrics. In particular, we focus on extremal metrics on the Fano surface $\mathbb{CP}^{2}\sharp2\overline{\mathbb{CP}}^{2}$. Work on this topic began around a decade ago with pioneering articles focusing on ideas of Simon Donaldson \cite{BuDo}, \cite{Do1}, \cite{Ke} (taking a more mathematical viewpoint) and, separately,  Matthew Headrick and Toby Wiseman \cite{Brau}, \cite{Dor}, \cite{HWK3}, \cite{HWDP2} (with a more physically motivated viewpoint).\\
\\
We propose a straightforward algorithm for the case of toric K\"ahler metrics which is easy to compute and overcomes some of the handicaps of existing techniques.  We remark that toric metrics were the subject of both \cite{BuDo} and \cite{Dor} but the toric condition is not integral to the algorithms that were employed there. The thrust of our method involves the minimization of a function of many variables and integration of a function of two variables over a polygon in the plane. The numerical techniques we use to achieve these objectives (conjugate gradient descent and Gaussian quadrature) are completely standard but nevertheless seem to achieve reasonable results. We hope that this article serves as a proof-of-concept and that it suggests more sophisticated numerical methods could yield even better accuracy. These ideas provide a good approximation to a distinguished extremal K\"ahler metric shown abstractly to exist by Chen, LeBrun and Weber \cite{CLW} (henceforth the extremal CLW metric). What distinguishes this particular extremal metric is that it is conformal to an Einstein metric. Throughout this paper we will refer to the K\"ahler metric as the extremal CLW metric and the conformally related Einstein metric as the Einstein CLW metric. As the extremal CLW metric does not have a closed form, finding numerical approximations is of considerable interest.\\ 
\\
A numerical approximation to the extremal CLW metric was given by Bunch and Donaldson in \cite{BuDo}. In general, Donaldson's algorithm considers metrics on polarised manifolds that are induced by embedding the manifold into a high-dimensional complex projective space equipped with the standard Fubini-Study metric; such metrics are called \emph{algebraic}.  Amongst these a distinguished algebraic metric known as a \emph{balanced metric} is found and then used as the starting point for a refined approximation. One significant drawback is that the procedure only works for metrics with rational K\"ahler classes. This is because it relies upon finding a Hermitian metric $h$ on a large tensor power of a holomorphic line bundle with the curvature of $h$ corresponding to the K\"ahler metric. The extremal CLW metric does not have a rational K\"ahler class (though the class is reasonably close to the canonical class and so an approximation can still be made).  Headrick and Wiseman's method is centred around being able to numerically simulate a parabolic flow, the K\"ahler-Ricci flow, that is expected to converge to the required metric. In this case the extremal metric is actually a K\"ahler-Einstein metric. The K\"ahler-Ricci flow can be framed as a second order parabolic PDE in a single function (the K\"ahler potential). For extremal metrics that are not K\"ahler-Einstein, the corresponding approach would be to utilize the Calabi flow. There are two immediate problems;  it is a fourth order PDE and so difficult to discretize, and it is not \emph{a priori} clear that one has convergence.\\ 
\\
Our method bypasses these problems. It is analogous, in some sense to the Donaldson method which involves minimizing functionals restricted to finite dimensional subspaces of metrics (the so-called algebraic metrics). We also minimize a functional restricted to a finite dimensional subspace of metrics, which we call the \emph{restricted symplectic metrics}. In fact, this method is suggested in the papers \cite{Dor} and \cite{HWDP2}.  Here the authors use the K\"ahler-Ricci flow to generate a representation of the metric. They then attempt to fit a restricted symplectic metric to this representation.  We do not need the first step and are able to search the space of restricted symplectic metrics directly. Our methods can also be seen as analogous to those employed in the recent work of Headrick and Nassar \cite{HN}. In this work the authors are concerned with numerically approximating Ricci-flat K\"ahler metrics on compact Calabi-Yau manifolds. They do this by considering a variety of different energy functionals (of varying order) and then minimising them on the space of algebraic metrics.  We also demonstrate a lower-order functional that seems to yield a good approximation to the extremal CLW metric when numerically minimised over the space of restricted symplectic metrics. We should remark that there is, as yet, no rigorous theoretical justification for the convergence of our method (in fact, only the Headrick-Wiseman method has what might be regarded as a satisfactory theory of convergence). In future work we hope to address some of the theoretical considerations surrounding the convergence of our algorithm. It would be particularly interesting to investigate whether the lower-order functional could be used as alternative to the Calabi Energy in the rigorous existence theory of the extremal CLW metric.  However, in this paper we will focus on the numerical approximations our method achieves.\\  
\\
An important application of numerical representations of canonical metrics has been to calculate various associated geometric invariants. Of particular interest is the calculation of the spectrum of certain natural differential operators such as the scalar Laplacian on functions or the Lichnerowicz Laplacian on symmetric 2-tensors.  One place where the spectrum of such operators appears is in the study of the Ricci flow
\begin{equation}\label{RF}
\frac{\partial g}{\partial t} = -2\Ric (g).
\end{equation}
Einstein manifolds (i.e. metrics $g$ such that  $\Ric (g) =\Lambda g$)  are fixed points of this flow; they only evolve via homothetic rescaling. A natural question to ask is whether a particular metric is stable as a fixed point (i.e. after a small perturbation does the Ricci flow return to the Einstein metric).  A result of Cao, Hamilton and Ilmanen \cite{CHI} (see \cite{CH} and \cite{HM14} for a proof) says that an Einstein metric is conformally unstable if the first non-zero eigenvalue of the ordinary Laplacian satisfies $\lambda_{1} < 2\Lambda$.  In \cite{HM14}, the authors investigated the torus-invariant spectrum, yielding the upper bound $\lambda_1 < 2.11\Lambda$.  A numerical proof of instability for the Einstein CLW metric was given in \cite{HHS}.  This was based on finding harmonic (1,1)-forms orthogonal to the K\"ahler form.  We recover these results using our approximations and give numerical evidence of an eigenvalue of the scalar Laplacian less than $2\Lambda$.\\
\\
Knowledge of the spectrum of the Lichnerowicz Laplacian of an Einstein metric can determine stability properties of related physical constructions \cite{GiHa}. Hence physicists are also interested in the stability of Einstein metrics. The numerical instability of a related non-K\"ahler, Hermitian, Einstein metric on $\mathbb{CP}^{2}\sharp\overline{\mathbb{CP}}^{2}$ due to Page \cite{Pa}, was first numerically demonstrated by Young \cite{You}. A theoretical proof of instability, not relying on numerical approximations, was given in \cite{HHS} and \cite{HM14}. Our numerical work suggests that the Einstein CLW metric is unstable in exactly the same manner as the Page metric.
\subsection{Computer Code}
The C++ code that this project uses to implement the conjugate gradient descent method of minimisation (Subsection \ref{CGD}) and the Matlab code used to implement the Levenberg-Marquardt algorithm (Subsection \ref{LM}) are available on both of the authors' webpages. \footnote{\text{http://www.buckingham.ac.uk/directory/dr-stuart-hall/}} \footnote{\text{http://mathfaculty.fullerton.edu/tmurphy/research.html}} The Matlab code contains fewer lines as we are using Matlab's `lsqnonlin' function in the optimisation toolbox to implement the Levenberg-Marquardt algorithm rather than coding our own version. Readers wishing to use the code should save the functions `CLWScal.m' and `CLWScalint.m' to a directory where Matlab can access them. Then, given an initial vector of inputs `$x0$', call the `lsqnonlin' function by typing \verb|x = lsqnonlin(@CLWScalint,x0)|
into the Matlab command line. This will return a vector `$x$' representing the optimised coefficients.\\
\\
\emph{Acknowledgements:} We are very grateful to Simon Donaldson, Matthew Headrick and Robert Haslhofer for providing comments on a draft version of this work. In particular, it was Matthew's suggestion to use the conjugate gradient method over the standard gradient descent. SH would like to thank Torben Kuseler for his assistance with running some of the algorithms. TM would like to thank McKenzie Wang for his support. We would like to thank the anonymous referees for their careful reading of the paper and for pointing out numerous corrections. Special thanks are due to the referee who drew our attention to the Levenberg-Marquardt algorithm. This has been a very useful step in helping along our related study of Ricci solitons and quasi-Einstein metrics on toric manifolds. This research was supported by a Dennison research grant.
\section{toric K\"ahler Metrics}
In this section we give a brief review of toric K\"ahler manifolds.  We refer the reader to the articles \cite{Abr1}, \cite{DonTor} for a comprehensive discussion of the theory. For our purposes a toric K\"ahler manifold will be a K\"ahler manifold $(M^{2n},\omega, J)$ that admits an effective action of the torus $\mathbb{T}^{n}$ that is simultaneously holomorphic and Hamiltonian. Crucially, there is a dense open subset $M^{o}\subset M$ on which this action is free. From the machinery developed in \cite{Abr1} and \cite{Gui} we obtain: 
\begin{itemize}
\item A compact convex polytope $P\subset \mathbb{R}^{n}$ and an identification $M^{o} \cong P^{o}\times \mathbb{T}^{n}$,
\item A finite set of  affine linear functions $l_{i}:\mathbb{R}^{n}\rightarrow \mathbb{R}$  such that the polytope $P$ is obtained as the intersection of the regions defined by $l_{i}(x)\geq 0$,
\item A smooth convex function $u:P^{o}\rightarrow \mathbb{R}$ such that in the coordinates on  $P^{o}\times \mathbb{T}^{n}$ the metric $g(\cdot,\cdot)=\omega(J\cdot,\cdot)$ has the form
$$g = u_{ij}dx_{i}dx_{j}+u^{ij}d\theta_{i} d\theta_{j}.$$
Here $u_{ij} = \frac{\partial^{2} u}{\partial x_{i} \partial x_{j}}$ and $u^{ij}$ is the ordinary matrix inverse.
\end{itemize}
The metric has coordinate singularities on the boundary $\partial P$ of $P$. It is known exactly how this has to occur. A result of Guillemin \cite{Gui}  and Abreu \cite{Abr1} shows that symplectic potential can be written as
\begin{equation}\label{Sympotdef}
u = \frac{1}{2}\sum_{i}l_{i}\log(l_{i})+F,
\end{equation}
where $F$ is a smooth function on $P$. The term $\frac{1}{2}\sum_{i}l_{i}\log(l_{i})$ is known as the canonical symplectic potential associated to $P$ and we will denote this by $u_{can}.$ Functions of the form (\ref{Sympotdef}) are said to satisfy the \emph{Guillemin boundary conditions}.
\subsection{Extremal metrics}
The K\"ahler metrics we are interested in approximating are known as extremal metrics.  They were introduced by Calabi \cite{Cal} as the critical points of the functional
$$\mathcal{C}(\omega) = \int_{M}S^{2} \ \frac{\omega^{n}}{n!}$$
where $S$ is the scalar curvature of $\omega$ and $\omega$ varies over all metrics in a fixed cohomology class.  The value of $\mathcal{C}(\omega)$ is referred to as the Calabi energy of the metric $\omega$. The Euler-Lagrange equations of the functional are equivalent to the requirement that $\nabla S$ is a holomorphic vector field, that is $\bar{\partial}\nabla S = 0.$
Viewed as a PDE in the K\"ahler potential, this is a highly non-linear sixth-order equation. This is one explanation as to why so little is known in general about these metrics. It is known that an extremal metric must be invariant under the maximal compact subgroup of the automorphism group. Hence an extremal metric must be torus invariant for any K\"ahler class that is toric. The scalar curvature of a toric K\"ahler metric is given by the following beautiful formula due to Abreu \cite{Abr2}
\begin{equation}
S=-u^{ij}_{ij}.
\end{equation}   
Another benefit of the toric setting is that the condition for a torus invariant function to have holomorphic gradient is that it is an affine linear function of the polytope coordinates. Putting things together we arrive at the following equation for an extremal metric on a toric K\"ahler manifold $(M^{2n}, \omega)$
\begin{equation} \label{Abreueqn}
u^{ij}_{ij} = \sum_{k=1}^{k=n} a_{k}x_{k}+b,
\end{equation}
where $a_{k}, b \in \mathbb{R}$. Equation (\ref{Abreueqn}) is usually referred to as \textit{Abreu's equation}.
A useful observation of Donaldson is that the constants $a_{k}$ and $b$ appearing in Abreu's equation can be determined from the polytope. In order to do this we need to define a measure $\sigma$ on the boundary of the polytope $\partial P$. This measure is just a multiple of the restriction of the Lesbegue measure on each edge scaled so that
$$|d\sigma\wedge dl_{r}| = dx$$
where $l_{r} $ is the affine linear function defining the $r^{th}$ edge and $dx$ is the standard Lesbegue measure on $\mathbb{R}^{n}$.
\begin{proposition}[c.f. Corollary 1 in \cite{DoInt}]
Let $u:P\rightarrow \mathbb{R}$ be a symplectic potential with Guillemin boundary conditions and with $S=-u^{ij}_{ij}$.  Then
\begin{equation}\label{DonIBP}
\int_{P}u^{ij}f_{ij} dx= \int_{\partial P}2fd\sigma-\int_{P}Sfdx,
\end{equation} 
for any $f\in C^{\infty}(P)$ that is continuous up to the boundary.
\end{proposition} 
The boundary term $\int_{\partial P}2f d\sigma$ is twice that which appears in the paper \cite{DoInt}. This is because the symplectic potential is dominated by $ \frac{1}{2}x \log (x) $ at the boundary of the polytope $P$. Donaldson uses a formulation of the theory where the singular behaviour is of the form $x \log (x)$. Equation (\ref{DonIBP}) shows that the quantity
$$\mathcal{L}_{S}(f) = \int_{\partial P}2fd\sigma-\int_{P}Sfdx$$ 
must vanish when $f$ is an affine-linear function. This places $n+1$ constraints on $S$ and so allows one to determine the affine-linear function $S$ exactly in the case of an extremal metric.
\subsection{Conformally K\"ahler, Einstein metrics}
One metric that we can apply our method to is distinguished amongst all extremal metrics on $\mathbb{CP}^{2}\sharp2\overline{\mathbb{CP}}^{2}$ in that it is conformal to an Einstein metric. The correspondence between Hermitian Einstein metrics on four-manifolds and extremal K\"ahler metrics was first noted by Derdzinski.
\begin{theorem}[Proposition 4 in \cite{Derd}]
Let $(M^{4},\omega,g)$ be a K\"ahler manifold of dimension 4, oriented in the natural way. Then the following conditions are equivalent:
\begin{enumerate}
\item The metric $S^{-2}g$, defined where the self-dual Weyl tensor $W^{+} \neq 0$, is an Einstein metric.
\item The metric $g$ has vanishing Bach tensor ($g$ is then said to be Bach-flat).
\end{enumerate}
Moreover, either of conditions (1) and (2) implies
\begin{equation}\label{confScal}
S^{3}+6S\Delta S-12|\nabla S|^{2}=\kappa,
\end{equation}
where $\kappa \in \mathbb{R}$ is the scalar curvature of the Einstein metric $g_{e}=S^{-2}g$.
\end{theorem}
We note that condition (1) automatically implies that the K\"ahler metric $g$ is extremal. The condition of being Bach-flat can also be interpreted as saying that the metric is a critical point of the Weyl curvature functional
$$\mathcal{W}(g) = \int_{M}|W(g)|^{2}dV_{g}$$ 
where $W(g)$ is the Weyl curvature tensor of $g$. On K\"ahler surfaces this functional is topologically equivalent to the Calabi energy. Hence a Bach-flat K\"ahler metric is one where the Calabi energy is extremised for nearby K\"ahler metrics. The fact that the extremal CLW metric is Bach-flat allows one to determine the parameter $a$ exactly.  Recently LeBrun was able to show that extremal  CLW metric globally minimises the Calabi energy for any K\"ahler metric on $\mathbb{CP}^{2}\sharp 2\overline{\mathbb{CP}}^{2}$ \cite{LB12}.\\
\\
As the scalar curvature of an extremal-toric metric can be explicitly calculated and the K\"ahler class of the extremal CLW metric is also known explicitly, one can compute the Einstein constant and hence the scalar curvature $\kappa=4\Lambda$ of the Einstein CLW metric appearing in Equation (\ref{confScal}).
\begin{lemma}[Lemma 2.4 in \cite{HM14}]
Let $(M^{4}, g_{e})$ be an Einstein metric satisfying $${\Ric (g_{e})=\Lambda g_{e}}.$$  Suppose further that $g_{e}=S^{-2}_{k}g_{k}$ for a K\"ahler metric $g_{k}$ with scalar curvature $S_{k}$. Then
$$\Lambda = \sqrt{\frac{96\pi^{2}\chi(M)+144\pi^{2}\tau(M)-\int_{M}S_{k}^{2} \ dV_{g_{k}}}{8Vol(g_{e})}},$$
where $\chi(M)$ is the Euler characteristic of $M$, $\tau(M)$ is the signature of $M$ and $Vol(g_{e})$ is the volume of the Einstein metric given by $\int_{M}S^{-4}_{k} \ dV_{g_{k}}$.
\end{lemma}
Using the description of the extremal CLW metric given in the next section, one can calculate $\kappa \approx 60.3456688$.
\subsection{The Chen-LeBrun-Weber metric}
The moment polytope here is a pentagon defined by the linear functions
$$ l_{1}(x) = 1+x_{1}, \ l_{2}(x) = 1+x_{2}, \ l_{3}(x)= a-1-x_{1}, \ l_{4}(x) = a-1-x_{2}$$ and $$ l_{5}(x) = a-1-x_{1}-x_{2}.$$
The constant $a$ determines the K\"ahler class of the metric.  If $a=2$ then the K\"ahler metric is in the class $c_{1}(M)$.  In \cite{LB95} LeBrun showed that the K\"ahler class of the extremal CLW metric, equal to $a\approx 1.9577128052$ to 10 d.p.  This is the value used in the numerics. In view of Equation (\ref{DonIBP}), the scalar curvature of an extremal metric in these classes is given by $${S=A(x_{1}+x_{2})+B},$$ where
 $$A=\frac{48(1-a^{3})}{a^{6}+6a^{5}+9a^{4}+4a^{3}-3a^{2}-6a+1}, $$
and
$$B= \frac{12(a^{5}+7a^{4}-2a^{3}+2a^{2}-5a+5)}{a^{6}+6a^{5}+9a^{4}+4a^{3}-3a^{2}-6a+1}.$$
We remark that these values are different from the ones calculated in \cite{HM14}.  This is for two reasons; firstly, there is a typographical error in the values given there (though all the calculations are performed with the correct values) and secondly, in this work we have translated the polytope by $(-1,-1)$.  This means that if the metric is in the class $c_{1}(M)$ (corresponding to $a=2$), the polytope agrees with the one used in \cite{HWDP2}.
\section{Algorithm}
The algorithm is centered around an expansion of the symplectic potential as
\begin{equation}\label{symppot}
u = u_{can}+F(x_{1},x_{2}) = u_{can}+\sum_{i,j}c_{ij}x_{1}^{i}x_{2}^{j}.
\end{equation}
We have not tried to find a rigorous justification for the symplectic potential of an extremal metric on a toric surface being real analytic in the polytope coordinates.  An argument in the special case of toric K\"ahler-Einstein metrics is given in \cite{Dor}. As our interest lies in the numerical results, we will suppress this technical point. In what follows we will describe the algorithm for approximating the extremal CLW metric though it is clear that one could perform the same procedure for any toric K\"ahler class. The extremal CLW metric is $\mathbb{Z}_{2}$ invariant (the action switches $x_{1}$ and $x_{2}$) so the function $F$ can be expanded as
$$F(x_{1},x_{2}) = c_{1}x_{1}x_{2}+c_{2}(x_{1}^{2}+x_{2}^{2})+c_{3}x_{1}x_{2}(x_{1}+x_{2})+c_{4}(x_{1}^{3}+x_{2}^{3})+...$$
Truncating the function $F$ by taking the first $n$ coefficients of the polynomial expansion means that the Calabi energy and the related integrals used, are functions of $n$ variables. The space of truncated representations is what we refer to as the space of \emph{restricted symplectic metrics}. Rather than trying to minimise the Calabi energy, we use the functional
\begin{equation}\label{MCE}
\mathcal{I}(\omega) = \int_{M}\big(S-S_{CLW}\big)^{2} \ \frac{\omega^{2}}{2},
\end{equation}
where $\omega$ ranges over the toric metrics in the same K\"ahler class as the extremal CLW metric and $S_{CLW}$ is the affine-linear function corresponding to the scalar curvature of the extremal CLW  metric (or the extremal metric in the K\"ahler class being considered). Restricting this to symplectic potentials of the form (\ref{symppot}) we obtain a function of $n$ variables, $\mathcal{I}_{n}$ which is given by
$$\mathcal{I}_{n}(c_{1},...,c_{n}) = 4\pi^{2}\int_{P}\big(S(c_{1},...,c_{n})-S_{CLW}\big)^{2} \ dx_{1}dx_{2},$$
where $S(c_{1},...,c_{n})$ is the scalar curvature of the metric given by the symplectic potential defined by the coefficients $c_{1},...,c_{n}$. We then proceed to minimise the functions $\mathcal{I}_{n} $ by the method of conjugate gradient descent.\\
\\
There are a variety of functionals that could be minimised in order to approximate the extremal CLW metric. We give a third-order functional in Section 4.2.  It is tempting to consider minimising the $L^{2}$-norm of the trace-free Ricci tensor of the conformal metric.  This would be a second-order functional and so it would appear easier to compute at first glance. The main disadvantage with this method is that the Ricci tensor has the same singular behaviour on the boundary of $P$ as the metric. Thus the calculation of the functional and its gradient becomes considerably more complicated. The scalar curvature is not singular on the boundary and so our algorithm does not require a particularly sophisticated integration scheme.  
\subsection{Conjugate Gradient Descent}\label{CGD}
This method is very widely used and we refer the reader to \cite{NR} for details.  We give an overview of the method here.  The rough idea is that one performs a gradient descent method, without repeating the search over directions that have already been tried. Consider first the problem of trying to minimise the quadratic function ${f: \mathbb{R}^{n} \rightarrow \mathbb{R}}$ given by
$$f(x) = c-b^{t}x +\frac{1}{2}x^{t}Ax,$$
where $b \in \mathbb{R}^{n}$ and $A$ is a symmetric positive-definite $n\times n$ matrix. It is easy to see that the minimum of $f$ is the solution to the equation $Ax=b$. The conjugate gradient algorithm in this case consists of making an initial guess $x_{0}$ which yields a residual vector ${r_{0} = b-Ax_{0}}$.  One then forms a sequence of vectors $r_{i}$ and $h_{i}$ defined by the recurrence
$$r_{i+1} = r_{i}-\alpha_{i}Ah_{i} \textrm{ and } h_{i+1} = r_{i+1}+\beta_{i}h_{i}$$
 where $p_{0}=r_{0}$ and the constants $\alpha_{i}, \beta_{i}$ are given by
$$\alpha_{i} = \frac{|r_{i}|^{2}}{h_{i}^{t}Ah_{i}} \textrm{ and } \beta_{i} = \frac{|r_{i+1}|^{2}}{|r_{i}|^{2}}.$$
Whilst carrying out this recurrence one updates the guess via
$$ x_{i+1} = x_{i}+\alpha_{i}h_{i}.$$
It is an exercise in linear algebra to show that this procedure finds the exact minimum of $f$ in at most $n$ steps.\\
\\
If  the smooth multivariable function ${F:\mathbb{R}^{n} \rightarrow \mathbb{R}}$ is well approximated by the quadratic function
$$Q_{F}(x) = F(a)+\nabla F (a)^{t}(x-a)+\frac{1}{2}(x-a)^{t}\nabla^{2} F(a) (x-a),$$
then it is natural to try to minimise $Q_{F}(x)$. The difficulty is that it would seem one needs to calculate the Hessian matrix $\nabla^{2} F(a)$ in order to carry out the algorithm described above. In fact, if at any point in the algorithm the residual $r_{i} =-\nabla F(P_{i})$ for some point $P_{i}$, then it is not difficult to show that the constant $\alpha_{i}$ is the value of $t$ that minimises the one-variable function
$$\tilde{F}(t)  = F(P_{i}+ t\cdot h_{i}).$$
Setting $P_{i+1} = P_{i}+\alpha_{i} h_{i}$ then one can also show that
$$r_{i+1} = -\nabla F (P_{i+1}).$$
Hence the algorithm for minimising the quadratic function $Q_{F}(x)$ can be implemented, without computing the Hessian of $F$, providing the gradient of $F$ can be computed and the one-variable functions $\tilde{F}(t)$ can be easily minimised.\\
\\
We implement the Polak-Ribere variant of the conjugate gradient method that is described on pages 518-519 of \cite{NR}.  We find the minimum of a one variable function by using Ridder's zero finding algorithm (described on page 453 of \cite{NR}) applied to the derivative. As the functions $\mathcal{I}_{n}$ are not exactly quadratic, we stop the algorithm after $4n$ steps and then restart at the latest guess.
\subsection{Integration - Gaussian Quadrature}
In order to calculate the various integrals we use the method of Gaussian quadrature (we refer the reader to \cite{NR} for more information).  For a one-dimensional integral (normalised so that the range is $[-1,1]$) the idea is to approximate the integral by taking a weighted sum of values
$$\int_{-1}^{1} f dx \approx \sum_{i=1}^{i=k}w_{i}f(x_{i}).$$
The points $x_{i}$ at which the function are sampled are known as the abscissa and the $w_{i}$ are referred to the weights. The points $x_{i}$ and weights $w_{i}$ are chosen so that if $f$ is a polynomial of degree $2k-1$ or less, then the sum will compute the integral exactly.
\\
To compute integrals over the pentagon $P$ we take the following splitting:
$$\int_{P}Fdx_{1}dx_{2} = \int_{-1}^{a-1}\int_{-1}^{1}Fdx_{1}dx_{2}+\int_{1}^{a-1}\int_{-1}^{a-1-x_{1}}Fdx_{2}dx_{1}.$$
This ensures that the all the one dimensional iterated integrals are smooth functions. We then approximate the iterated one-dimensional integrals using the Gaussian quadrature method. In order to check the accuracy of the methods we compute the volume of the Einstein CLW metric 
$$Vol(g_{CLW}) = 4\pi^{2}\int_{P}S_{CLW}^{-4}dx_{1}dx_{2} \approx 0.583421245.$$ 
Both the Gaussian quadrature with 10 and 20 points agree with this value up to 9 decimal places. Hence we use the 10 point Gaussian quadrature to compute the integrals appearing in the algorithm.
 
\section{Results}
The algorithm was implemented in C++ using the value $a = 1.9577128052$. For each degree,  the previous coefficients were used as the initial guess with the value $0$ being entered where no previous value had been calculated. The conjugate gradient method was carried out twenty-five times or until no change in the Calabi energy was noticeable to nine decimal places.  Various measures of the accuracy of the numerical approximations were calculated.  The most obvious is essentially the value of the functional we are trying to minimise. The appropriate measure is the volume-normalised $L^{2}$-norm of the difference between the scalar curvature of the approximation $S$ and the affine-linear function $S_{CLW}$ representing the scalar curvature of the extremal CLW metric $g_{CLW}$, i.e.
$$ \left(\frac{1}{Vol(g)}\int_{M}\big(S-S_{CLW}\big)^{2} \ dV_{g_{k}}\right)^{1/2}=\left(\frac{2}{a^{2}+2a-1}\int_{P}\big(S-S_{CLW}\big)^{2} \ dx_{1}dx_{2}\right)^{1/2}.$$
As we know the extremal  CLW metric to be conformal to an Einstein metric (and we know \emph{a priori} the Einstein constant)  we calculate some measure of this discrepancy. The measure we use is the volume-normalised $L^{2}$-norm of
$$ \kappa-S_{CLW}^{3}-6S_{CLW}\Delta S_{CLW}+12|\nabla S_{CLW}|^{2} .$$
This is a measure of how far the conformal metric $g_{e} = S^{-2}g$ is from having constant scalar curvature and obviously if the metric is the extremal CLW metric then that above quantity vanishes (cf. Equation (\ref{confScal})).\\
\\
The final measure of accuracy we use is to compare certain integrals of the gradient of the scalar curvature to values that the authors were able to calculate in a closed form (i.e. a rational function of the parameter $a$) in \cite{HM14}. We have the following result allowing the calculation of $\|\nabla S_{CLW}^{p}\|_{L^{2}}$ (here the gradient and the $L^{2}$-norm are with respect to the Einstein CLW metric). 
\begin{proposition}[Proposition 2.3 in \cite{HM14}]
Let $(M^{4},g_{k})$ be a Riemannian manifold with everywhere nonzero scalar curvature $S_{k}$ and let $g_{e}=S^{-2}_{k}g_{k}$.  Let $\kappa$ be the scalar curvature of the metric $g_{e}$.  Then for $p\neq 1/2$ we have
$$\int_{M}|\nabla_{e}S_{k}^{p}|^{2}dV_{g_{e}} = \frac{p^{2}}{6(2p-1)}\int_{M}(S_{k}^{4}-\kappa S_{k})S_{k}^{2p-5}dV_{g_{k}}.$$
\end{proposition}
\begin{cor}
Using $a=1.9577128052$ and performing all calculations with the Einstein CLW  metric, we obtain
$$\|\nabla S_{CLW}\|^{2}_{L^{2}} \approx 4.9689665$$
and 
$$\|\nabla S_{CLW}^{-1}\|^{2}_{L^{2}} \approx 0.020806979. $$
\end{cor}
Given this result, we compute the values of $\|\nabla_{u} S_{CLW}\|^{2}_{L^{2}}$ and $\|\nabla_{u} S^{-1}_{CLW}\|^{2}_{L^{2}}$ with respect to the metric $S^{-2}_{CLW}g$, where $g$ is the restricted symplectic metric coming from our approximation to the extremal CLW metric.  More explicitly, if
$${S_{CLW} =A(x_{1}+x_{2})+B}$$ and the symplectic potential of $g$ is denoted by $u$, then
$$\|\nabla _{u}S_{CLW}\|^{2}_{L^{2}} =  4\pi^{2}\int_{P}A^{2}(u^{11}+2u^{12}+u^{22})S_{CLW}^{-2}dx_{1}dx_{2},$$
and
$$ \|\nabla_{u} S_{CLW}^{-1}\|^{2}_{L^{2}} = 4\pi^{2}\int_{P} A^{2}(u^{11}+2u^{12}+u^{22})S_{CLW}^{-6}dx_{1}dx_{2}. $$

The results of the numerical search are included in Table \ref{T1}. The \textbf{Deg} column refers to the degree of truncated polynomial. The column \textbf{$L^{2}$-error} shows the volume normalised $L^{2}$ difference between the scalar curvature and the affine-linear function $S_{CLW}$. The column  \textbf{Max} is the maximum pointwise difference between the scalar curvature and the affine-linear function $S_{CLW}$ and the column \textbf{Min} is the minimum difference. The column $\beta$ is the volume-normalised $L^{2}$ difference in the conformal scalar curvatures.  The columns $\|\nabla_{u} S_{CLW}\|^{2}_{L^{2}}$ and $\|\nabla _{u}S_{CLW}^{-1}\|^{2}_{L^{2}}$ are the $L^{2}$-norms of the gradient of $S_{CLW}$ and $S_{CLW}^{-1}$ computed with the approximate metric.

\begin{table}[!ht]
\centering
\caption{Errors and related quantities for scalar curvature minimisation procedure}
\begin{tabular}{|c|c|c|c|c|c|c|}
\hline
\textbf{Deg} &  \textbf{$L^{2}$-error} & \textbf{Max} & \textbf{Min} & \textbf{$\beta$} & $\|\nabla_{u} S_{CLW}\|^{2}_{L^{2}}$ & $\|\nabla_{u} S_{CLW}^{-1}\|^{2}_{L^{2}}$\\
\hline
\hline
2 &  0.48 & 3.64 & -1.08 & 1.9 & 4.751605 & 0.0196836\\
\hline
3 &  0.25 &2.25 & -0.33 & 0.77 & 4.931697 & 0.0205739\\
\hline
4 &  0.13 & 1.42 & -0.81 & 0.35 & 4.956189 & 0.0207408\\
\hline
5  & 0.066 & 0.90 & -0.13 & 0.15 & 4.964912 & 0.0207880 \\
\hline
6  & 0.035 & 0.57 & -0.34 & 0.070 & 4.967752 & 0.0208026 \\
\hline
7  & 0.019 & 0.35 & -0.054 & 0.033 & 4.968622 & 0.0208056 \\
\hline
8  & 0.010 & 0.21 & -0.13 & 0.016 & 4.968868 & 0.0208065\\
\hline
9  & 0.0052 & 0.12 & -0.020 & 0.0075 & 4.968939 & 0.0208068\\
\hline
10 & 0.0026 & 0.072 & -0.045 & 0.0038 & 4.968959 & 0.020806952 \\
\hline
11  & 0.0014 & 0.043 & -0.0090 & 0.0019 & 4.968963 & 0.020806964  \\
\hline
12 & $7.3 \times 10^{-4}$ & 0.024 & -0.021 & $9.2\times 10^{-4}$ & 4.9689659  & 0.020806976 \\
\hline
13 & $3.9 \times 10^{-4} $ & 0.014 & -0.0026  & $5.0\times 10^{-4}$ & 4.9689651  & 0.020806971 \\
\hline
14 & $2.1 \times 10^{-4}$ & 0.0086 & -0.0025 & $2.6 \times 10^{-4}$ & 4.9689659 & 0.020806976  \\
\hline
15 & $1.2\times10^{-4}$ & 0.0048 & -0.0011 & $1.4\times10^{-4}$ & 4.9689664 & 0.020806977 \\
\hline
16 & $7.4\times10^{-5}$ & 0.003 & $-6.8\times10^{-4}$ & $9.5\times10^{-5}$ & 4.9689652  & 0.020806972  \\
\hline
17 & $4.9\times 10^{-5}$ & 0.0021 & $-5.4\times 10^{-4}$ & $5.7\times 10^{-5}$ & 4.9689657 & 0.020806975 \\
\hline
18 & $4.4\times 10^{-5}$ & 0.0019 & $-4.9\times 10^{-4}$ & $5.3\times 10^{-5}$ & 4.9689650 & 0.020806973 \\
\hline
19 & $3.5\times 10^{-5}$ & 0.0014 & $-3.9\times 10^{-4}$ & $5.0\times 10^{-5}$ & 4.9689643 & 0.020806970  \\
\hline
20 & $2.5\times 10^{-5}$ & 0.0010 &  $-3.6\times 10^{-4}$ & $3.6\times 10^{-5}$ & 4.9689666 & 0.020806978 \\
\hline
\end{tabular}

\label{T1}
\end{table}

\vspace{10pt}
As the table shows, the approximations seem to be converging to the extremal CLW metric. The calculations of the gradients of $S_{CLW}$ for the higher degree approximations seem to agree with the exact values to 6 significant figures. It is reasonable to expect that the eigenvalue calculations in section 5 are also accurate to this level. We remark that the degree 20 approximation can be achieved in a couple of hours on a standard desktop computer. A combination of more sophisticated hardware and better algorithms could probably yield faster convergence.\\
\\
As with the numerical approximation to the K\"ahler-Ricci soliton on $\mathbb{CP}^{2}\sharp 2\overline{\mathbb{CP}}^{2}$ in \cite{HWDP2}, we give the explicit quartic approximation to the extremal CLW metric.
$$ u = u_{can} -0.09962x_{1}x_{2}-0.1333(x_{1}^{2}+x_{2}^{2})$$
$$-0.04195x_{1}x_{2}(x_{1}+x_{2}) - 0.03139(x_{1}^{3}+x_{2}^{3})$$
$$-0.01471x_{1}^{2}x_{2}^{2}-0.01119x_{1}x_{2}(x_{1}^{2}+x_{2}^{2})-0.007613(x_{1}^{4}+x_{2}^{4}).$$
\subsection{Comparison to the Bunch-Donaldson approximation}
In \cite{BuDo} the authors give an approximation to an extremal K\"ahler metric in the first Chern class of ${M=\mathbb{CP}^{2}\sharp 2\overline{\mathbb{CP}}^{2}}$ which corresponds to taking the parameter $a=2$.  The method is to notice that such a metric can always be be described as the curvature of a Hermitian metric $h$ on the anticanonical line bundle $K^{-1}$. A metric on the line bundle induces the usual $L^{2}$ metric on the space of global holomorphic sections of powers of $K^{-1}$, the finite dimensional vector space $H^{0}(M,K^{-p})$, and so induces a Fubini-Study metric on $\mathbb{P}(H^{0}(M,K^{-p}))$. For large enough powers of $p$ there is an embedding of $M$ into $\mathbb{P}(H^{0}(M,K^{-p}))$ by orthonormal sections. The pullback of the Fubini-Study metric gives (after dividing by $p$) a K\"ahler metric in $c_{1}(M)$. This process is iterated until a fixed point is reached.  The fixed point is known as a balanced metric. For large values of $p$ the balanced metrics approximate the extremal metric to order $O(p^{-1})$ in any $C^{k}$-norm and so Bunch and Donaldson then run a refined approximation algorithm starting at the balanced metric.  Essentially they use an elegant Newton-Raphson type algorithm to minimise the functional 
$$\mathcal{I}(\omega) = \int_{M}\big( S-S_{CLW}\big)^{2} \ \frac{\omega^{2}}{2}$$ 
restricted to the space of algebraic metrics which are those those coming from embedding $M$ into $\mathbb{P}(H^{0}(M,K^{-p}))$.\\
\\
The representation of an algebraic K\"ahler metric can be realised as a Hermitian $N_{p}\times N_{p}$ matrix where $N_{p}=\dim H^{0}(M,K^{-p})$. In this case $N_{p} = O(p^{2})$ and for example, when $p=20$, $N_{p} = 190$. As the metrics considered are torus invariant this means that the matrices are diagonal. Furthermore, as the extremal metric in this class is $\mathbb{Z}_{2}$-invariant, the metric can be represented by $95$ real coefficients. The representation of the extremal metrics by way of the polynomial expansion of the symplectic potential requires $\lfloor\frac{n^{2}+4n-4}{4}\rfloor$ real coefficients for a degree $n$ representation. So for example when the potential is a degree $20$ polynomial this involves a representation of the metric using $119$ real coefficients. It appears our approximation is not as efficient as the Bunch-Donaldson approximation, as when the degree is 20, our normalised $L^{2}$-error is $2.5\times10^{-5}$ whilst they are able to achieve $1.02 \times 10^{-6} $ when $p=20$.
\subsection{A third-order functional}
We also investigated an approximation technique based around minimising another functional other than the modified Calabi energy. This is somewhat analogous to the investigations carried out by Headrick and Nassar \cite{HN} where they use a variety of functionals to find approximations to Ricci-flat K\"ahler metrics on Calabi-Yau manifolds. We use the functional derived from the identity (\ref{confScal})
$$\mathcal{J}(\omega) = \int_{M}\left( \kappa-S_{CLW}^{3}-6S_{CLW}\Delta S_{CLW}+12|\nabla S_{CLW}|^{2}  \right)^{2} dV_{g}.$$
As the terms involving the metric, $\Delta S_{CLW}$ and $|\nabla S_{CLW}|^{2}$, only need one derivative of the metric, $\mathcal{J}$ is a third-order functional.  This makes the gradient significantly faster to compute. We carry out the conjugate gradient descent method as with the modified Calabi energy. The results of this minimisation are contained in Table \ref{T2}.

\begin{table}[!ht]
\centering
\caption{Errors and related quantities for conformal scalar curvature minimisation procedure}
\begin{tabular}{|c|c|c|c|c|c|c|c|}
\hline
\textbf{Deg} &  \textbf{$L^{2}$-error} & \textbf{Max} & \textbf{Min} & \textbf{$\beta$} & $\|\nabla_{u} S_{CLW}\|^{2}_{L^{2}}$ & $\|\nabla_{u} S_{CLW
}^{-1}\|^{2}_{L^{2}}$\\
\hline
\hline
2 &   0.5165  & 4.15 & -1.25 & 1.852 & 4.813795 & 0.0199199\\
\hline
3 &   0.2828 & 2.78 & -0.43 & 0.6772 & 4.916611 & 0.0204998 \\
\hline
4 &  0.1475  & 1.78 & -0.65 & 0.2820 & 4.957289 & 0.0207248 \\
\hline
5  & 0.08018 & 1.18 & -0.15 & 0.1200 & 4.964739 & 0.0207798 \\
\hline
6  & 0.04545 & 0.75  & -0.28 & 0.05412 & 4.967980 & 0.0207996 \\
\hline
7  & 0.02428 &  0.50  & -0.064  &  0.02493  & 4.968575 & 0.0208044 \\
\hline
8  & 0.01434 & 0.31  & -0.11 &  0.01155  & 4.968869 & 0.0208062 \\
\hline
9  &  0.007127  & 0.19 & -0.026 & 0.005370 & 4.968923 & 0.0208067 \\
\hline
10  & 0.004872 & 0.12 & -0.026 & 0.002853 & 4.968950 & 0.0208069 \\
\hline
\end{tabular}

\label{T2}
\end{table}

As the table shows, the results are in line with the approximations obtained by minimising the modified Calabi energy. The advantage of this method is that the algorithm seems to take far fewer steps to converge and each of the steps involves fewer calculations.

\subsection{Nonlinear least squares methods}\label{LM}

One referee of the paper suggested that it is useful to consider the approximations of the integrals $\mathcal{I}_{n}$ given by Gaussian quadrature as a sum of squares functional. The 20 point procedure can be thought of as choosing 800 points in the polytope (as we split the integral into two parts) and then evaluating a nonlinear sum of squares. Hence
$$\mathcal{I}_{n}(c_{1},...,c_{n}) \approx \sum_{i=1}^{i=800}\tilde{w}_{i}(S(c_{1},...,c_{n})-S_{CLW})^{2}(p_{i},q_{i}),$$
where $\tilde{w_{i}}$ is the appropriate weight and $(p_{i},q_{i}) \in P^{\circ} $. Such least-squares problems are particularly amenable to a method of optimisation known as the Levenberg-Marquardt algorithm. Matlab uses this method in the `lsqnonlin' function which appears as part of the optimisation toolbox. We implemented this algorithm in Matlab. The function tolerance and step-size tolerance were set to be $1\times 10^{-13}$ (i.e. the algorithm terminates if the change in the residual or in the approximate value of $c$ has absolute value less than $1\times 10^{-13}$), the maximum number of function evaluations was set at at 6000.  All other parameters were left at the default settings.\\
\\
The results follow exactly the pattern of Table \ref{T1} hence we can be very confident that the restricted symplectic metrics produced by the algorithm are converging to the extremal CLW metric. There are enormous benefits to using Matlab over the C++ routine. The main one is that one only needs to code in the scalar curvature function for a particular $(c_{1},...,c_{n})$ and not the actual implementation of the Levenberg-Marquardt algorithm which is contained in the `lsqnonlin' function. Had we been aware of this method when we began our project we would probably not written any C++ code to implement the conjugate gradient descent algorithm. However both algorithms seem to achieve the same results and so we present both methods here.\\
\\
Further details of the Levenberg-Marquardt are available in \cite{NR}. The authors have also used this approach to investigate other canonical metrics (gradient Ricci solitons and quasi-Einstein metrics) on the manifolds $\mathbb{CP}^{2}\sharp \overline{\mathbb{CP}}^{2}$ and $\mathbb{CP}^{2}\sharp 2 \overline{\mathbb{CP}}^{2}$. Details of this investigation and of the algorithm will appear in \cite{HM15}.

\section{Applications}
One use of the numerical approximations to canonical metrics is the calculation of various geometric information.  A geometric invariant of particular interest is the first non-zero eigenvalue of the scalar Laplacian. This is of particular importance for Einstein metrics as it can determine whether or not the Einstein metric is linearly stable as a fixed point of the Ricci flow.\\
\\
In \cite{HM14}, the fact that Einstein CLW metric is conformal to an extremal metric was used to give an explicit upper bound for the first non-zero eigenvalue of the scalar Laplacian of the Einstein metric, $\lambda^{CLW}_1$.  Expanding functions as powers of the scalar curvature $S_{CLW}$ the authors were able to show $\lambda^{CLW}_{1}\leq 2.107\Lambda$. We note that the space of functions on $\mathbb{CP}^{2}\sharp 2\overline{\mathbb{CP}}^{2}$ can be decomposed into the $\mathbb{Z}_{2}$-invariant and $\mathbb{Z}_{2}$-anti-invariant functions. We denote the first space by $\mathcal{F}^{+}$ and the second by $\mathcal{F}^{-}$ and note that the Laplacian preserves the decomposition $C^{\infty}(M)=\mathcal{F}^{+}\oplus\mathcal{F}^{-}$. The space $\mathcal{F}^{+}$ can be further decomposed into $\mathcal{F}^{+} = \mathbb{R}\oplus\mathcal{F}^{+}_{0}$ where $\mathcal{F}^{+}_{0}$ is the space of functions with vanishing integral. We consider minimising the Rayleigh quotient $\frac{\|\nabla f\|^{2}_{L^{2}}}{\|f\|^{2}_{L^{2}}}$ for cubic functions $f \in \mathcal{F}^{+}_{0}$ and $f\in \mathcal{F}^{-}$. More precisely, we consider minimising the Rayleigh quotient over the functions
$$\phi^{+} = c +(x_1+x_2)+a_{1}x_{1}x_{2} +a_{2}(x_{1}^{2}+x_{2}^{2})+a_{3}x_{1}x_{2}(x_{1}+x_{2})+a_{4}(x_{1}^{3}+x_{2}^{3})$$
and 
$$ \phi^{-} = (x_{1}-x_{2})+a_{1}(x_{1}^{2}-x_{2}^{2})+a_{2}x_{1}x_{2}(x_{1}-x_{2})+a_{3}(x_{1}^{3}-x_{2}^{3}).$$
The constant $c$ is chosen to ensure that $\phi^{+} \in \mathcal{F}_{0}^{+}$.
Using the degree 20 approximation we minimise the Rayleigh quotient over $(a_{1},a_{2},a_{3},a_{4})$ for $\phi^{+}$ and $(a_{1},a_{2},a_{3})$ for $\phi^{-}$. This yields

\begin{table}[!ht]
\centering
\caption{Values of coefficients for approximate eigenfunctions}
\begin{tabular}{|c|c|c|}
\hline
$a_{i}$ & $\phi^{+}$ & $\phi^{-}$ \\
\hline
\hline
$a_{1}$ & 0.7241 & 0.2894 \\
\hline
$a_{2}$ & 0.3141 & 0.1133 \\
\hline
$a_{3}$ & 0.1829 & 0.0774 \\
\hline
$a_{4}$ & 0.0790 & - \\
\hline
\end{tabular}

\end{table}

This results in values for the Rayleigh quotients of
$$ \frac{\|\nabla\phi^{+}\|^{2}_{L^{2}}}{\|\phi^{+}\|^{2}_{L^{2}}} \approx 2.0940\Lambda \textrm{ and } \frac{\|\nabla\phi^{-}\|^{2}_{L^{2}}}{\|\phi^{-}\|^{2}_{L^{2}}} \approx 1.9481\Lambda.$$
We note that the value of $2.0940\Lambda$ is very close to the value $2.0965\Lambda$ found by the authors in \cite{HM14} by considering polynomials in $(x_{1}+x_{2})$ where the integrals could be evaluated in a closed form. This is further evidence that the approximation we have given is close to the true extremal CLW  metric.
\subsection{Linear Stability}
It is clear that Einstein metrics can be considered as fixed points of the Ricci flow (\ref{RF}) as they evolve only by homothetic scaling. It is a natural question to ask whether they are attracting or repelling as fixed points. One way of determining this is to use Perelman's monotone quantity $\nu(g)$ \cite{Per}. This quantity is monotonically increasing under the Ricci flow and constant only if the metric is a gradient Ricci soliton which is a metric solving the equation
\begin{equation}\label{GRS}
\Ric (g)+ \nabla ^{2} f = \lambda g.
\end{equation}
The notion of Ricci soliton generalises that of an Einstein metric (the Einstein condition being recovered by setting $f$ constant). At an Einstein metric $g_{e}$ one can compute
$$ \frac{d^{2}}{dt^{2}}\nu(g_{e}+th) |_{t=0},$$
where $h\in s^{2}(TM)$.  If this quantity is positive then the Einstein metric $g_{e}$ is unstable (linearly unstable) as a small perturbation of the metric in the direction of $h$ will never flow back to $g_{e}$. The stability of an Einstein metric is related to the spectrum of the Lichnerowicz Laplacian, ${\Delta_{L}:s^{2}(TM)\rightarrow s^{2}(TM)}$, where 
$$\Delta_{L}h = \Delta h-2Rm(h,\cdot)+\Ric \cdot h+h\cdot \Ric.$$
This is not surprising as the Lichnerowicz Laplacian is essentially the linearization of the Ricci tensor viewed as a differential operator on symmetric 2-tensors. The following theorem makes this precise.
\begin{theorem}[Theorem 2.1 in \cite{CHI}, Theorem 1.2 in \cite{CZ}]
Let $(M,g_{e})$ be an Einstein manifold with $\Ric(g_{e})=\Lambda g_{e}$.  Then if $\Delta_{L} \geq 2\Lambda$ (i.e. the smallest non-zero eigenvalue of $\Delta_{L}$ is greater than $2\Lambda$), $g_{e}$ is linearly stable as a shrinking Ricci soliton.
\end{theorem}
In order to compute the spectrum of $\Delta_{L}$ at an Einstein metric one can use the fact that there are a number of subbundles that it preserves. In particular, if one considers conformal perturbations of the metric then the stability criterion can be phrased in terms of the spectrum of the ordinary Laplacian.
\begin{theorem}[Proposition 2.6 in \cite{CHI}, Theorem 1.1 in \cite{CH}]
Let $(M^{n},g)$ be an Einstein metric satisfying $\Ric (g)=\Lambda g$ and let $\lambda_{1}$ be the first non-zero eigenvalue of the scalar Laplacian. Then if
$$\frac{n\Lambda}{n-1}<\lambda_{1}<2\Lambda,$$
$g$ is linearly unstable as a shrinking Ricci soliton.
\end{theorem}
Hence, from the above results, we have numerical evidence for the following:
\begin{conj}
The Einstein CLW metric is linearly unstable and can be destabilised by a conformal perturbation. 
\end{conj}
Using the fact that the Einstein CLW metric is conformal to an extremal K\"ahler metric and $h^{(1,1)}(M)=3$ one can consider perturbations in the direction of harmonic $(1,1)$-forms.  Here there is the following theorem.
\begin{theorem}[Theorem 6.3 and Lemma 6.4 in \cite{HHS}] Let $(M^{4},g_{e})$ be an Einstein metric with scalar curvature $\kappa$, conformal to a K\"ahler metric $g_{k}$ by $g_{e}=S^{-2}_{k}g_{k}$. Suppose that the K\"ahler structure on the manifold has $h^{(1,1)}(M)>1$. Then if $\Delta_{k}S^{2}_{k}<\frac{\kappa}{2}$, $(M,g_{e})$ is linearly unstable.
\end{theorem}
In fact the authors prove that if  $\Delta_{k}S^{2}_{k}<\frac{5}{16}\kappa$ then the Ricci-flat cone that can be constructed over the Einstein CLW metric is unstable. Using Equation (\ref{confScal}) and the identity $\Delta_{k}S_{k}^{2} = 2S_{k}\Delta_{k} S_{k}  -2|\nabla S_{k}|^{2}$ we have
$$\Delta_{k}S^{2}_{k} = \frac{\kappa}{3}+2|\nabla {S}_{k}|^{2}-\frac{S_{k}^{3}}{3}.$$
Using the degree 20 approximation we compute that
$$\Delta S^{2}_{CLW} \leq 10.33 < \frac{5}{16}\kappa \approx 18.86.$$
Hence we recover the evidence of \cite{HHS} that the Einstein CLW metric is linearly unstable as a shrinking soliton and that the Ricci-flat cone over it is also unstable.


\begin{thebibliography}{999999}
\bibitem{Abr2} M. Abreu, \textit{K\"ahler geometry of toric varieties and extremal metrics}, Internat. J. Math., 9(6), 641--651 (1998).
\bibitem{Abr1} M. Abreu, \textit{K\"ahler geometry of toric manifolds in symplectic coordinates}, Fields Inst. Commun., 35, AMS, Providence, Rhode Island, 1--24 (2003).
\bibitem{Brau} V. Braun, T. Brelidze, M. R. Douglas, B. A. Ovrut, \textit{Eigenvalues and eigenfunctions of the scalar Laplace operator on Calabi-Yau manifolds}, J. High Energy Physics, Art No. 120, pp 1--57 (2008). 
\bibitem{BuDo} R. S. Bunch, S. K. Donaldson, \textit{Numerical approximations to extremal metrics on toric surfaces}, Handbook of geometric analysis, Adv. Lect. Math. 7, no. 1, 1–-28 (2008).
\bibitem{Cal} E. Calabi, \textit{Extremal K\"ahler metrics}, Ann. of Math. Stud., 102, Princeton Univ. Press, Princeton N. J., 259--290 (1982).
\bibitem{CHI} H.-D. Cao, R. Hamilton, T. Ilmanen, \textit{Gaussian density and stability for some Ricci solitons}, arXiv:0404.165 (2004).
\bibitem{CH}  H.-D. Cao, C. He, \textit{Linear stability of Perelman's $\nu$-entropy on symmetric spaces of compact type}, J. Reine Angew. Math (to appear) (2013).
\bibitem{CZ}  H.-D. Cao, M. Zhu, \textit{On second variation of Perelman's shrinker entropy}, Math. Ann., 353(3), 747--763 (2012).
\bibitem{CLW} X. X. Chen, C. LeBrun, B. Weber, \textit{On conformally K\"ahler, Einstein manifolds}, J. Amer. Math. Soc., 21(4), 1137--1168 (2008).
\bibitem{Derd} A. Derdzinski, \textit{Self-dual K\"ahler manifolds and Einstein manifolds of dimension four}, Compos. Math., 49, 405--433 (1983).
\bibitem{DoInt} S. K. Donaldson, \textit{Interior estimates for solutions of Abreu's equation}, Collect. Math., 56(2), 103--142 (2005).
\bibitem{Do1}  S. K. Donaldson, \textit{Some numerical results in complex differential geometry}, Pure Appl. Math. Q. 5, no. 2, 571–-618 (2009).
\bibitem{DonTor} S. K. Donaldson, \textit{K\"ahler geometry on toric manifolds and some other manifolds with large symmetry}, Handbook of Geometric Analysis, Adv. Lect. Math., 7(1), 29--75 (2008).
\bibitem{Dor} C. Doran, M. Headrick, C. Herzog, J. Kantor, T. Wiseman, \textit{Numerical K\"ahler-Einstein on the third del Pezzo}, Commun. Math. Phys., 282(2), 357--393 (2008).
\bibitem{NR} B.P. Flannery, W. H. Press, S. A. Teukolsky, W.T. Vetterling, \textit{Numerical recipes: the art of scientific computing - third edition}, Cambridge University Press, (2007). 
\bibitem{GiHa} G. Gibbons, S. Hartnoll, \textit{Graviational instability in higher dimensions}, Phys. Rev. D, 66, 24--64 (2002). 
\bibitem{Gui} V. Guillemin, \textit{K\"ahler structures on toric varieties}, J. Diff. Geom., 40(2), 285--309 (1994).
\bibitem{HHS} S. J. Hall, R. Haslhofer, M. Siepmann, \textit{The stability inequality for Ricci-flat cones}, J. Geom. Anal., 24(1), 472--494 (2014). 
\bibitem{HM14} S. J. Hall, T. Murphy, \textit{On the spectrum of the Page and Chen-LeBrun-Weber metrics}, Ann. Glob. Anal. Geom., 46(1), 87--101 (2014). 
\bibitem{HM15} S. J. Hall, T. Murphy, \textit{Approximating Ricci solitons and quasi-Einstein metrics on toric surfaces}, (in preparation) (2015).
\bibitem{HN} M. Headrick, A. Nassar, \textit{Energy functionals and Calabi-Yau metrics}, Adv. Theor. Math. Phys. 17, 867--902 (2013). 
\bibitem{HWK3} M. Headrick, T. Wiseman, \textit{Numerical Ricci-flat metrics on $K3$}, Class. Quantum Gravity, 22(3), 4931--4960 (2005). 
\bibitem{HWDP2} M. Headrick, T. Wiseman, \textit{Numerical K\"ahler-Ricci soliton on the second del Pezzo}, arXiv:0706.2329v1, (2007).
 \bibitem{Ke} J. Keller, \textit{Ricci iterations on K\"ahler classes}, J. Inst. Math. Jussieu, 8(4), 743--768, (2009).
\bibitem{LB95}  C. LeBrun, \textit{Einstein metrics on complex surfaces}, In: Pedersen, H., Andersen, J. Dupont, J., Swann A. (eds.) Geometry and Physics (Aarhus, 1995). Lecture Notes in Pure and Applied Mathematics, vol. 184, pp. 167–176. Dekker, New York (1997).
\bibitem{LB12} C. LeBrun, \textit{On Einstein Hermitian 4-manifolds}, J. Diff. Geom., 90, 277-302 (2012).
\bibitem{Pa} D. Page, \textit{A compact rotating gravitational instanton},  Phys. Lett., 79B, 235--238 (1979).
\bibitem{Per} G. Perelman, \textit{The entropy formula for the Ricci flow and its geometric applications},  arXiv:math/0211159v1 (2002).
\bibitem{You} R. E. Young, \textit{Semiclassical instability of gravity with positive cosmological constant}, Phys. Rev. D 28(10), 2436–-2438 (1983).
\end{thebibliography}
\end{document}